\documentstyle[12pt]{article}
\newcommand{\const}{\mathop{\rm const}\limits}

\newcommand{\Var}{\mathop{\rm Var}\limits}

\newcommand{\Law}{\mathop{\rm Law}\limits}

\newcommand{\Sub}{\mathop{\rm Sub}\limits}

\newcommand{\SSub}{\mathop{\rm SSub}\limits}

\newcommand{\ASSub}{\mathop{\rm ASSub}\limits}

\newcommand{\Bin}{\mathop{\rm Bin}\limits}

\begin{document}

 \begin{center}

{\bf  NON-UNIFORM NON-ASYMPTOTICAL SHARP  \\

\vspace{3mm}

 ESTIMATE FOR THE RATE OF CONVERGENCE  \\

\vspace{3mm}

FOR BERNSTEIN'S POLYNOMIAL APPROXIMATION, } \\

\vspace{4mm}

{\bf with bilateral constant evaluation}\\

\vspace{4mm}

{\sc Eugene Ostrovsky, Leonid Sirota}\\

\vspace{3mm}

 Bar-Ilan University,  59200, Ramat Gan, ISRAEL; \\
 e-mail: eugostrovsky@list.ru \\
 e-mail: sirota3@bezeqint.net \\

\vspace{4mm}

     {\sc Abstract.}

\end{center}

 \ We  derive the non-asymptotical non-uniform sharp error estimation for Bernstein's approximation of
continuous function based on the modern probabilistic apparatus.\par
 \ We investigate also the convergence of derivative of these polynomials and we will consider briefly also the
multivariate case. \par

\vspace{4mm}

\par

\vspace{3mm}

\section{ Introduction. Notations. Statement of problem.}

\vspace{3mm}

  \ Let $  f = f(x), \ x \in [0,1] $ be continuous: $  f \in C[0,1], $ numerical function,  $  n = 2,3, \ldots, $

$$
B_n[f](x) = \sum_{m=0}^n  {n \choose m } \ f \left(\frac{m}{n} \right) \ x^m \ (1 - x)^{n-m} \eqno(1.0)
$$
be its Bernstein's polynomial of degree $  n. $ \par
 \ It is reasonable to define formally

$$
B_{\infty}[f](x) = f(x), \ f(\cdot) \in C[0,1]. \eqno(1.0a)
$$

 \ Denote also as ordinary

$$
\omega[f](\delta) = \sup_{x,y \in [0,1], |x - y| \le \delta} |f(x) - f(y)| \ - \eqno(1.1)
$$
 module (modulus) of continuity of the function $  f = f(x) $ of a first order,

$$
\Delta_n[f](x) = | \ B_n[f](x) - f(x) \ |  \eqno(1.2)
$$
be an {\it non-uniform} error of the Bernstein's approximation of a $ n^{th} $ order  for the function $  f. $ \par

 \ The problem  of error estimation for Bernstein's approximation of continuous function goes back to the classical work
belonging to  S.N.Bernstein \cite{Bernstein1} \ (1912). He proved in the famous article \  \cite{Bernstein1} \ that as $  n \to \infty $
the sequence $ \ B_n[f](x) \ $ of a  Bernstein's polynomial of degree $  n $ converges uniformly as $ n \to \infty $
to the initial function $  f(x): $

$$
\lim_{n \to \infty} \max_{x \in [0,1]} \left| \  B_n[f](x) - f(x) \ \right| = \lim_{n \to \infty} \max_{x \in [0,1]} \Delta_n[f](x) = 0.
$$

\ There exists a huge number of publications about the estimation of a value $ \Delta_n[f](x), $
 see for example  \cite{Achiezer1}, \cite{Bojanic1}, \cite{Chlodovsky1}, \cite{DeVore1}, \cite{Li2}, \cite{Lorentz1}, \cite{Mat'e},
 \cite{Palatanea1}, \cite{Prestin1}, \cite{Telyakovskii1}, \cite{Tikhonov1}, \cite{Timan1} and many others; see also
 reference therein. \par

\vspace{4mm}

{\bf  Our purpose in this short report is obtaining a sharp up to multiplicative constant non-uniform and non-asymptotical
universal estimate for the error of the Bernstein's approximate for the arbitrary initial continuous function in the terms
of its module of continuity of the first order.} \par

 \vspace{3mm}

 \ We intend to use too modern probabilistic apparatus: theory of Grand Lebesgue Spaces of random variables,
in particular, theory of Strictly Subgaussian random variables; and introduce and use some generalization of this
notion: Almost Strictly Subgaussian random variables. \par

\vspace{4mm}

 \ Note that the non-uniform estimate are "much" better in the considered case than uniform ones, see
 \cite{Palatanea1}, \cite{Palatanea2}, \cite{Lorentz1} . \par

\vspace{3mm}

\section{ Auxiliary apparatus: subgaussian and strictly subgaussian random variables.}

\vspace{3mm}

 \ Let $  \{\Omega, B, {\bf P}  \} $ be some sufficiently rich probability space with expectation $ {\bf E.} $ \\

 \vspace{2mm}
 {\bf Definition  2.1.}\\
 \vspace{2mm}

  \  We say that the {\it centered:} $ {\bf E} \xi = 0 $ numerical random variable (r.v.)
 $ \xi $ is subgaussian, or equally, belongs to the space $ \Sub(\Omega), $
 if there exists some non-negative constant $ \tau \ge 0 $ such that

$$
\forall \lambda \in R  \ \Rightarrow
{\bf E} \exp(\lambda \xi) \le \exp[ \lambda^2 \ \tau^2 /2]. \eqno(2.1).
$$

\vspace{3mm}

 \ The minimal value $ \tau $ satisfying (2.1) for all the values $ \ \lambda \in R \ $  is called a  {\it subgaussian  norm}
of the variable $ \xi, $ write

 $$
 ||\xi||\Sub = \inf \{ \tau, \ \tau > 0: \ \forall \lambda \in R \ \Rightarrow {\bf E}\exp(\lambda \xi) \le \exp(\lambda^2 \ \tau^2/2) \}.
 $$

 Evidently,

$$
||\xi||\Sub = \sup_{\lambda \ne 0} \left[ \sqrt{ 2 \ln {\bf E}  \exp (  \lambda \xi)  }/|\lambda| \right].  \eqno(2.2)
$$

 \ This important notion was introduced by  J.P.Kahane \cite{Kahane1}; V.V.Buldygin and Yu.V.Kozachenko in \cite{Buldygin1} proved
that the set $ \Sub(\Omega) $  relative the norm $  ||\cdot|| $ is complete Banach space which is isomorphic to subspace
consisting only from the centered variables of Orlicz's space over $ (\Omega, B,P)  $ with $ N \ -$ Orlicz-Young function
 $ N(u) = \exp(u^2) - 1; $  see also \cite{Kozatchenko1}.  \par

\vspace{3mm}

 \  The detail investigation of this class or random variables with very interest  applications into the theory of random fields
  reader may found in the book \cite{Buldygin3};  we reproduce here some main facts from this monograph.\par

\vspace{3mm}

{\bf Definition 2.1.} The tail function $  T_{\eta}(u), \ u > 0   $ for the numerical valued random variable $  \eta  $
 is as usually defined as follows

$$
T_{\eta}(u) \stackrel{def}{=} \max [{\bf P}(\eta > u), \ {\bf P}(\xi < -u) ], \ u > 0.
$$

\vspace{3mm}

  \ If for instance $ ||\xi||\Sub = \tau \in (0,\infty),  $ then

 $$
T_{\xi}(u) =
 \max [{\bf P}(\xi > u),  {\bf P}(\xi < -u)  ] \le \exp(- u^2/(2 \tau^2)  ), \ u \ge 0; \eqno(2.3)
 $$
 and  the last inequality is in general case non-improvable.  It is sufficient for this to consider the case when
 the r.v. $  \xi  $ has the centered Gaussian non-degenerate distribution.\par

  \ Conversely, if  $ {\bf E} \xi = 0 $ and if  for some positive finite constant $  K  $

 $$
  T_{\xi}(u) = \max [{\bf P}(\xi > u),  {\bf P}(\xi < -u)  ] \le \exp(- u^2/K^2  ), \ x \ge 0,
 $$
 then $ \xi \in \Sub(\Omega) $ and $ ||\xi||\Sub < 4 K. $ \par

 \ The subgaussian norm in the subspace of the centered r.v. is equivalent to the following Grand Lebesgue Space (GLS)
 norm:

 $$
|||\xi||| := \sup_{s \ge 1} \left[ \frac{|\xi|_s}{\sqrt{s}} \right], \hspace{6mm} |\xi|_s =  \left[ {\bf E} |\xi|^s \right]^{1/s}.
 $$

 \ More detail investigation of these spaces see in the monograph \cite{Ostrovsky1}, chapter 1.  \par

  \ Denote in the sequel for brevity  for any r.v.  $ \eta $

 $$
 \sigma^2(\eta) = \sigma^2 = \Var{\eta} = {\bf E}\eta^2 - ( {\bf E} \eta)^2.
 $$

\vspace{3mm}

{\bf Definition 2.2.} (See \cite{Buldygin3}, chapter 1; \ \cite{Ostrovsky102}.) \ The subgaussian r.v. $ \xi $ is said to be
 {\it Strictly Subgaussian,} write $  \xi \in \SSub, $ iff

$$
\forall \lambda \in R \ \Rightarrow {\bf E} e^{\lambda \xi} \le e^{\lambda^2 \sigma^2(\xi)/2}, \eqno(2.4)
$$
or equally

$$
||\xi ||\Sub \le \sigma(\xi) = ||\xi||L_2(\Omega). \eqno(2.4a)
$$

 Recall that always $ ||\xi ||\Sub \ge \sigma(\xi) = ||\xi||L_2(\Omega),  $ so that

$$
\xi \in \SSub(\Omega) \Leftrightarrow  {\bf E}\xi = 0, \ ||\xi ||\Sub = \sigma(\xi) = ||\xi||L_2(\Omega).  \eqno(2.4b)
$$

 \ Many examples of strictly subgaussian distributions may be found in the book of  V.V.Buldygin and Yu.V.Kozatchenko \cite{Buldygin3},
chapter 1. For instance, arbitrary mean zero Gaussian distributed r.v. is strictly subgaussian, including the case when
this r.v. is equal to zero  a.e.; the symmetric Rademacher's r.v. $ \rho $ with distribution
$ {\bf P}(\rho = 1) = {\bf P}(\rho = -1) = 1/2  $  belongs to the set  $ \SSub(\Omega). $
The random variable $  \eta $ which has an uniform distribution on the symmetrical interval $ (-b,b), \ b = \const \in (0,\infty) $ is
Strictly Subgaussian. \par

 \ Consider also for instance following the authors \cite{Buldygin3}  the r.v. $ \zeta $  with the following density:

 $$
 f_{\zeta}(x) =   \frac{\alpha+1}{2 \alpha} \  \left(1 - |x|^{\alpha} \right) \ I(|x| \le 1), \
 \alpha = \const \ge 0, \eqno(2.5)
 $$
 where $ I(A) = I(A,x) = 1, \ x \in A; \ I(A) = I(A,x) = 0, \ x \notin A $ is indicator function;
then $ \zeta \in \SSub(\Omega). $ \par

 \ This example is interesting because the kurtosis of the r.v. $ \zeta $ is zero if $ \alpha = \sqrt{10} - 3. $

\vspace{3mm}

 The convenience of these notions is following. Let $  \{ \xi(i ) \}, \ i = 1,2,\ldots,n  $ be (centered) independent
subgaussian r.v. Denote

$$
S(n) = \sum_{i=1}^n \xi(i), \hspace{6mm} \Sigma^2(n) = \sum_{i=1}^n (||\xi(i)||\Sub)^2. \eqno(2.6)
$$

 Then $ ||S(n)||\Sub \le \Sigma(n) $ and following

$$
 \max( {\bf P}(S(n)/\Sigma(n) > x ),  {\bf P}(S(n)/\Sigma(n) < -x )) \le e^{ -x^2/2  }, \ x \ge 0, \eqno(2.7)
$$
 the tail or concentrations inequalities. \par
 \  If in addition $ \xi(i) $ are identical distributed and $ \beta:= ||\xi(1)||\Sub \in (0,\infty), $ then

 $$
 \sup_n ||S(n)/\sqrt{n}||\Sub = \beta
 $$
  and

$$
\sup_n \max( {\bf P}(S(n)/(\beta \sqrt{n} ) > x ),( {\bf P}(S(n)/(\beta \sqrt{n} ) < - x ) \le e^{ -x^2/2  }, \ x \ge 0, \eqno(2.8)
$$

 If in addition the r.v. $ \xi(i) $  are strictly subgaussian, the estimate (2.8) may be reinforced by lower estimate  used the classical CLT:

$$
\sup_n {\bf P}(S(n)/(\beta \sqrt{n} ) > x ) \ge \lim_{n \to \infty} {\bf P}(S(n)/(\beta \sqrt{n} ) > x ) =
$$

$$
( 2 \pi )^{-1/2} \int_x^{\infty} e^{-y^2/2} \ dy \ge C  \ x^{-1} e^{-x^2/2}, \ x \ge 1. \eqno(2.9)
$$

\vspace{4mm}

{\bf Definition 2.3.}The non-degenerate  centered  random variable $ \nu  $ with variance  $  \sigma^2,
0 < \sigma < \infty $ is said to be {\it  almost strictly subgaussian,} briefly $ \nu \in \  {\bf \ASSub,} \  $ if
for all the real values $  \lambda \in R $

$$
{\bf E} \cosh(\lambda \nu/\sigma) \le \exp \left(\lambda^2/2 \right). \eqno(2.10)
$$

\vspace{3mm}

 \ {\bf Proposition 2.1.} \par

\vspace{3mm}

 \  It follows immediately from the definition (2.10) the following tail estimate for these variables
by means of Tchebychev-Chernov inequality

$$
T_{\nu}(u)  \le 2 \exp \left(-u^2/2 \sigma^2 \right), \ u \ge 0. \eqno(2.11)
$$

\vspace{3mm}

 \ {\bf Proof.}  \ Obviously, if the r.v $  \nu $ is almost strictly subgaussian: $ \nu \in {\bf \ASSub}  $ then
$  c \cdot \nu \in {\bf \ASSub}.   \ c = \const. $ \par
 \ Let now $ \nu \in {\bf \ASSub} $ with $ \ \sigma = 1. $ We have by means of Tchebychev's  inequality for the positive
values $  \lambda $ and $  u: $

$$
T_{\nu}(u) \le  \frac{e^{\lambda^2/2}}{\cosh (\lambda u)} = \frac{e^{\lambda^2/2}}{0.5 ( e^{\lambda u} + e^{-\lambda u} )} \le
$$

$$
2 \exp(\lambda^2/2 - \lambda u) = 2 \exp(-u^2/2),
$$
if we choose $ \lambda = u. $ \par

 \ Note that analogous approach appears at first in an article \cite{Ostrovsky103}, where was applied to the investigation
of the Central Limit Theorem in Banach space. \par

\vspace{4mm}

 \ Denote also

$$
 \theta = \theta(p) = \sqrt{ p \cdot (1-p)}, \ 0 < p < 1. \eqno(2.12)
$$

\vspace{3mm}

 \ Let us now formulate and prove the main result of this section.  Denote by $  \mu = \mu_n = \mu_{n,p} $ the r.v. having
Bernoulli (Binomial) distribution: $  \Law \mu = \Bin(n,p), \ n = 1,2,\ldots; \ 0 < p < 1. $ Recall that

$$
{\bf P} (\mu = m) =  {n \choose m } \ p^m \ (1 - p)^{n-m}, \ m = 0,1,2,\ldots,n
$$
and

$$
{\bf E} \mu = np, \ \Var \mu = n \ p \ (1-p) = n \ \theta(p).
$$

\vspace{3mm}

 \ {\bf Remark 2.1.} The non-centered random variable $  \zeta $ may be named Almost Strictly Subgaussian iff
it has a first moment and the centered r.v. $ \zeta - {\bf E} \zeta $ satisfies the definition 2.3. \par

\vspace{3mm}

{\bf  Theorem 2.1 } (in our definitions and notations). \par

\vspace{3mm}

 {\it  The random variables } $  \mu - np, \ 0 < p < 1, \ n = 1,2,\ldots $  {\it are almost strictly subgaussian
uniformly for all the values }  $  n. $ \par

\vspace{3mm}

{\bf Proof.}  The inequality

$$
{\bf E}(\mu - n p)^{2m} \le n^{-m} \frac{(2m)!}{2^m \ m!} \ \theta^m(p), \ m = 0,1,2, \ldots \eqno(2.13)
$$
is proved in the book of G.G.Lorentz \cite{Lorentz1}, page 14. We get denoting the centered and normed r.v.

$$
\eta :=   \frac{\mu - n p}{\sqrt{n p (1 - p)}}
$$
for all the real values  $  \lambda \in R: $

$$
{\bf E} \cosh (\lambda \eta) \le \sum_{m=0}^{\infty} \frac{\lambda^{2m}}{(2m)!} \ {\bf E} \eta^{2m} \le
$$

$$
\sum_{m=0}^{\infty} \frac{\lambda^{2m}}{(2m)!}  \frac{(2m)!}{2^m \ m!}  =
\sum_{m=0}^{\infty} \frac{\lambda^{2m}}{2^m \ m!} =  e^{ \lambda^2/2 }. \eqno(2.14)
$$

\vspace{4mm}

 \ {\bf Remark 2.2.}  The inequality

$$
\ln {\bf E} e^{\lambda \mu_1} \ln \left[(1 - p) e^{-\lambda p} + p e^{\lambda (1-p) }  \right] \le p(1-p) \cdot \frac{\lambda^2}{2},
\eqno(2.15)
$$
where $ 1/2 \le p \le 1  $  and $  \lambda \ge 0 $
is proved at first in an article  of D.Berend  and  A.Kontorovich \cite{Berend1}, lemma 5, page 4;
see also some applications in \cite{Schlemm1} and a preliminary results in \cite{Buldygin2},\cite{Buldygin3}, \cite{Kearns1}. \par

 \ Perhaps,  this form of the probabilistic  inequality was known for the specialists in the approximation theory, see e.g.
\cite{Gonska1}, \cite{DeVore1}, chapter 10, \cite{Palatanea1}, \cite{Palatanea2}. \par

\vspace{3mm}

\section{ Main result. Exactness.}

\vspace{3mm}

  \ Let us introduce the following (sublinear) operator (more precisely, the sequence of operators)

$$
J_n[f](x) \stackrel{def}{=} \int_0^{\infty} \omega[f] \left( \frac{z \cdot \theta(x)}{\sqrt{n}} \right) \
 z  \ \exp \left( - \frac{z^2}{2} \right) \ dz.
$$

 \ Denote also

$$
W = W(x) \stackrel{def}{=}
 \sup_{ \const \ne f \in C[0,1]} \sup_{n = 1,2,\ldots}  \left[ \frac{\Delta_n[f](x)}{J_n[f](x) } \right].   \eqno(3.0)
$$

\vspace{3mm}

{\bf Theorem 3.1.}  {\it  We assert that for all the values  } $  x \in (0,1)  $

$$
\frac{1}{\pi} \le W = W(x) \le 2.  \eqno(3.1)
$$

\vspace{3mm}

{\bf Proof. Upper bound.} \\

 \ Let $ x \in (0,1);  $ both the degenerate cases $ x = 0 $ or $  x = 1 $ are trivial. \par

  \ The  expression for $  B_n[f](x)  $ may be represented likewise in the initial work of   S.N.Bernstein  \cite{Bernstein1}
 as follows

$$
B_n[f](x)  = {\bf E} f(\mu/n), \ \eqno(3.2)
$$
where the r.v. $  \mu $ has a binomial (Bernoulli)  distribution with parameters $ {\bf E} \mu = n x, \ \Var \mu = n x (1-x). $
 The r.v. $  \mu - nx $ is equal to the sum

$$
\mu - n x = \sum_{j=}^n \xi_x(j), \ \Leftrightarrow \mu/n -  x = \sum_{j=1}^n \xi_x(j)/n, \eqno(3.3)
$$
where  the r.v. $ \xi_x $  has a binary distribution

$$
{\bf P} (\xi_x = 1 - x) = x = 1 - {\bf P} (\xi_x = - x), \ x \in (0,1)
$$
and the random variables $  \xi_x(j)  $ are independent copies of the r.v. $ \xi_x. $ \par

 \ Put

 $$
 \zeta = \zeta_{n,x}  = \sum_{j=1}^n \xi_x(j)/\sqrt{n}, \eqno(3.4)
 $$
then $  \mu/n = x + \zeta_{n,x}/\sqrt{n}.  $ As we knew, the linear combination of strictly subgaussian independent random
variables are also strictly subgaussian and we conclude therefore on the basis of theorem 2.1

$$
T_{\zeta_{n,x}}(u) \le 2 \ \exp \left( -u^2/2 \theta^2(x) \right). \eqno(3.5)
$$

 \ We can and will suppose without loss of generality that the r.v. $ \zeta_{n,x} $ has a symmetrical distribution
such that

$$
T_{\zeta_{n,x}}(u) = 2 \ \exp \left( -u^2/2 \theta^2(x) \right). \eqno(3.5a)
$$

 \ More precisely, let $  \tilde{\zeta} = \tilde{\zeta}_{n,x} $ be a symmetrically distributed r.v. such that

$$
T_{\tilde{\zeta}_{n,x}}(u) = 2 \ \exp \left( -u^2/2 \theta^2(x) \right), \eqno(3.5b)
$$
then for any positive monotonically increasing function $ G = G(z), \ z \ge 0  $

$$
{\bf E} G(\zeta_{n,x}) \le  {\bf E} G(\tilde{\zeta}_{n,x}).
$$

 \ We estimate consequently

$$
\Delta_n[f](x) = \left| \ B_n[f](x) - f(x) \ \right| =
$$

$$
 \left| {\bf E}(f(x + \zeta/\sqrt{n}) - f(x)) \right| \le
{\bf E} \omega[f] \left( \frac{\zeta_{n,x}}{\sqrt{n}}  \right),
$$
and we deduce after integration by parts relaying on the inequality (3.5)

$$
\Delta_n[f](x) \le 2 \ \int_0^{\infty} \omega[f] \left( \frac{y}{\sqrt{n}}  \right) \ \frac{y}{\theta^2(x)} \
\exp \left( - y^2 /2 \theta^2(x)  \right) \ dy =
$$

$$
2 \ \int_0^{\infty} \omega[f] \left( \frac{z \theta(x)}{\sqrt{n}}  \right) \ z
\exp \left( - z^2 /2 \right) dz = 2 \ J_n[f](x). \eqno(3.6)
$$

 \ Thus, we proved the upper bound: $  W \le 2.  $ \par

\vspace{4mm}

{\bf Remark 3.1} Note that $ \theta(x) \le 1/2,  $ therefore

$$
\Delta_n[f](x) \le 2 \int_0^{\infty} \omega[f] \left( \frac{y}{2 \sqrt{n}}  \right) \ y \
\exp \left( - y^2 \right) \ dy. \eqno(3.6a)
$$

 \ Thus, we obtained easily  the uniform Bernstein's error estimate through the uniform one. \par

\vspace{4mm}

{\bf Remark 3.2} It follows immediately from the estimation (3.6a) and on the basis of Lebesgue dominated
convergence theorem that

$$
\lim_{n \to \infty} \max_{x \in [0,1]} \Delta_n[f](x) = 0,
$$
i.e. the classical Bernstein's result. \par

\vspace{4mm}

{\bf Examples 3.1.}  \\

\vspace{2mm}

 \ Let the function $  f = f(x) $ be H\"olderian:

$$
\exists \alpha = \const \in (0,1], \ H_{\alpha} = \const < \infty \ \Rightarrow
$$

$$
|f(x_1) -  f(x_2)| \le H_{\alpha} \ |x_1 - x_2|^{\alpha},  \ x_1, x_2 \in [0,1]; \eqno(3.7)
$$
or equally

$$
\omega[f](\delta) \le H_{\alpha} \ \delta^{\alpha}, \hspace{4mm} H_{\alpha} < \infty.  \eqno(3.7a)
$$

 \ It follows immediately from theorem 3.1

$$
\Delta_n[f](x) \le 2 \ H_{\alpha} \ \left[  \frac{2x(1-x)}{n}  \right]^{\alpha/2} \  \Gamma \left( \frac{\alpha}{2}  \right) . \eqno(3.8)
$$

 \ If in addition $ \alpha = 1  $ (Lipshitz condition) and if we denote as ordinary $  L = H_1 \in (0,\infty),  $ then

$$
\Delta_n[f](x) \le 2 \ L \ \left[  \frac{2 \pi x(1-x)}{n}  \right]^{1/2}. \eqno(3.8a)
$$

 \ See also \cite{Mat'e}, \ \cite{Prestin1}. \par

\vspace{3mm}

{\bf Lower bound.}\\

\vspace{2mm}

  \ Ranko Bojanic in the article \cite{Bojanic1}, see also \cite{Bojanic2}, \cite{Tikhonov1}, \cite{Telyakovskii1} considered the
family of a trial functions of a form

$$
g_t(x) = |t-x|, \ t,x \in [0,1],
$$
which if Lipshitzian relative the variable $  x, $ (as well as relative the variable $  t), $ with constant $  L = 1, $
and proved  in fact that as $  n \to \infty  $

$$
\Delta_n[g_t(\cdot)] = \left[  \frac{2 x(1-x)}{\pi n}  \right]^{1/2} + 0\left( \frac{1}{n} \right). \eqno(3.9)
$$

 \  We conclude therefore that for all the values $ x $ from the (open) interval $  x \in (0,1)  $

$$
W(x) \ge \lim_{n \to \infty}
  \left[ \frac{\Delta_n[g_t](x)}{J_n[g_t](x) } \right]  = 1/\pi.   \eqno(3.10)
$$

 \ This completes the proof  of theorem 3.1. \par

\vspace{3mm}

\section{ Convergence of derivative for Bernstein's polynomials.}

\vspace{3mm}

  \ Let now in this section the aforementioned function $  f $ be continuous differentiable: $  f \in C^1[0,1].  $
It is known that the derivatives  $ B'[f](x)  $ converges uniformly to one for the function $  f'(x), $ see
\cite{Floater1}, \cite{Veretennikov1}.\par

 \ We intend to obtain here the refined estimation also non-asymptotical and non-uniform for the approximation

$$
\Delta^1_n[f](x) \stackrel{def}{=} \max_x \left| B'_n[f](x)  - f'(x) \right|. \eqno(4.1.)
$$

\vspace{3mm}

{\bf Theorem 4.1.}

$$
\Delta^1_n[f](x) \le \frac{3}{2} \omega[f'] \left( \frac{1}{n} \right) + 2 \ J_{n-1} [f'](x) =
$$

$$
 \frac{3}{2} \omega[f'] \left( \frac{1}{n} \right) +
 2 \ \int_0^{\infty} \omega[f'] \left( \frac{z \cdot \theta(x)}{\sqrt{n - 1}} \right) \
 z  \ \exp \left( - \frac{z^2}{2} \right) \ dz, \ n \ge 2. \eqno(4.2)
$$

\vspace{3mm}

{\bf Proof.} \par

\vspace{3mm}

 \ Let $  n \ge 2. $  The expression for $   B'[f](x) $ is given, for example, in the article
\cite{Veretennikov1}:

$$
B'_n[f](x) = \sum_{j=0}^{n-1} n \ V_{1/n}[f] \left(\frac{j}{n} \right) \ {n-1 \choose j } \ x^j (1-x)^{n-1-j}, \eqno(4.3)
$$
where

$$
V_{1/n}[f] \left(\frac{j}{n} \right) =   f \left(\frac{j+1}{n} \right) - f\left(\frac{j}{n} \right)
$$
is a difference operator for the function $  f $  with step $  1/n  $  at the point $ j/n.  $ \par

 \ Note that

$$
 \left| \ n \ V_{1/n}[f] - f' \left( \frac{j}{n} \right)  \ \right| \le \frac{1}{2} \ \omega[f'] \left( \frac{1}{n} \right). \eqno(4.4)
$$

  \ Further, since $  n \ge 2 $

$$
\left| f' \left( \frac{j}{n} \right)  - f' \left( \frac{j}{n-1} \right)  \right| \le \omega[f'] \left( \frac{1}{n}  \right),
$$
therefore (see theorem 3.1)

$$
\Delta^1_n[f](x) \le \frac{3}{2} \omega[f'] \left( \frac{1}{n} \right) + \left| \ B_{n - 1}[f'](x) - f'(x) \ \right| \le
$$

$$
\frac{3}{2} \omega[f'] \left( \frac{1}{n} \right) + 2 \ J_{n-1} [f'](x), \eqno(4.5)
$$
 Q.E.D. \par

\vspace{3mm}

{\bf Remark 4.1.}

 \ As for the lower estimate for Bernstein derivative estimate, consider the trial functions of the form

$$
G_t = G_t(x) := \int_0^x g_t(y) dy; \eqno(4.6)
$$
then for all the values $ t \in [0,1] $

$$
\Delta^1_n[G_t](x) \ge \int_0^{\infty} \omega[G_t'(x)] \left( \frac{z \cdot \theta(x)}{\sqrt{n}} \right)
\ z \ \exp(-z^2/2) \ dz. \eqno(4.7)
$$

\vspace{3mm}

\section{Multivariate Bernstein's polynomials. }

\vspace{3mm}

\begin{center}

{\bf The speed of convergence.} \par

\end{center}

\vspace{3mm}

 \ Let now $  f = f(x,y), \ (x,y) \in [0,1]^2  $ be continuous numerical function defined on the  square $ [0,1]^2, \
n_1, n_2 \ge 2 $ be integer numbers.

 \ One can define the modulus of continuity of the multivariate function $ f = f(x,y) $ as follows

$$
\omega[f](\delta_1, \delta_2) \stackrel{def}{=}
\sup_{|x_1 - y_1|\le \delta_1} \sup_{|x_2 - y_2|\le \delta_2} |f(x_1, y_1) - f(x_2, y_2)|, \ \delta_{1,2} \in [0,1]. \eqno(5.0)
$$

 \ The multivariate (more exactly, bivariate)
 Bernstein's polynomial $  B_{n_1,n_2}[f](x,y)   $ of an order $  (n_1, n_2)  $ one can define by a formula

$$
B_{n_1,n_2}[f](x,y) \stackrel{def}{=}
$$

$$
\sum_{k_1=0}^{n_1} \ \sum_{k_2=0}^{n_2} {n_1 \choose k_1 } {n_2 \choose k_2 }
f \left( \frac{k_1}{n_1}, \ \frac{k_2}{n_2}\right)\ x^{k_1} (1 - x)^{n_1 - k_1} \ y^{k_2} (1 - y)^{n_2 - k_2}, \eqno(5.1)
$$
see, e.g. \cite{Lorentz1}, \cite{Veretennikov1}. \par

  \ We define formally by the analogy with the one-dimensional case

$$
B_{n_1,\infty}[f](x,y) \stackrel{def}{=}
 \sum_{m=0}^{n_1}  {n_1 \choose m } \ f \left(\frac{m}{n_1}, y \right) \ x^m \ (1 - x)^{n-m} \eqno(5.1a)
$$
and analogously

$$
B_{\infty, n_2}[f](x,y) \stackrel{def}{=}
 \sum_{k=0}^{n_2}  {n_2 \choose k } \ f \left(x, \frac{k}{n_2} \right) \ y^k \ (1 - y)^{n_2-k} \eqno(5.1b).
$$

  \ Finally,

 $$
 B_{\infty, \infty} [f](x,y)\stackrel{def}{=} f(x,y), \ f (\cdot, \cdot) \in C \left([0,1]^2  \right).
 $$

 \ Define also the two-dimensional non-uniform error of Bernstein approximation

$$
\Delta_{n_1, n_2}[f](x,y):= | \ B_{n_1, n_2}[f](x,y) - f(x,y) \ |.  \eqno(5.2)
$$

\vspace{3mm}

 \ Introduce analogously to the one-dimensional case the following  operator
(more precisely, the sequence of operators) $ J_{n_1, n_2}[f](x,y) \stackrel{def}{=} $

$$
\int_0^{\infty} \int_0^{\infty} \omega[f] \left( \frac{z_1 \cdot \theta(x)}{\sqrt{n_1}}, \  \frac{z_2 \cdot \theta(y)}{\sqrt{n_2}} \right)
 z_1 \ z_2  \ \exp \left( - \frac{z_1^2}{2} -  \frac{z_2^2}{2} \right) \ dz_1 \ dz_2. \eqno(5.3)
$$

\vspace{4mm}

{\bf Theorem 5.1.}  {\it  We propose as before that for all the values  } $ x,y \in [0,1]  $ {\it  and } $  d = 2 $

$$
\pi^{-1} \le W_2 = W_2(x,y) \stackrel{def}{=}
 \sup_{ \const \ne f \in C[0,1]^2} \sup_{n_{1,2} = 1,2,\ldots}  \left[ \frac{\Delta_{n_1,n_2}[f](x,y)}{J_{n_1, n_2}[f](x,y) } \right]  \le 2^d.   \eqno(5.4)
$$

\vspace{3mm}

{\bf Proof.} \par

\vspace{3mm}

 \ Our reasoning are likewise ones in the section 3. Namely, let $ 0 < x,y < 1; $ we have

$$
B_{n_1, n_2}[f](x,y) = {\bf E} f\left( \frac{\mu_1}{n_1}, \frac{\mu_2}{n_2} \right), \eqno(5.5)
$$
where the {\it independent} random variables $  \mu_1  $  and $  \mu_2 $  have the integer binomial (Bernoulli)
distribution with the parameters correspondingly

$$
{\bf E } \mu_1 = n_1 x, \ {\bf E } \mu_2 = n_2 y, \ \Var \mu_1 = n_1 x(1-x), \  \Var \mu_2 = n_2 y(1-y). \eqno(5.6)
$$

 \ Each such the variable may be represented as above

$$
\frac{\mu_1}{n_1} = x + \frac{\zeta_{1,x}}{\sqrt{n_1}}, \hspace{5mm} \frac{\mu_2}{n_2} = y + \frac{\zeta_{2,y}}{\sqrt{n_2}}, \eqno(5.7)
$$
where both the r.v. $ \zeta_{1,x} $  and $  \zeta_{2,y} $ are independent and are as above almost strictly subgaussian. \par

 \ We can and will suppose without loss of generality that the r.v. $ \zeta_{1,x}, \ \zeta_{2,y} $ have a symmetrical distribution
such that

$$
T_{\zeta_{1,x}}(u_1) = \exp \left( -u_1^2/2 \theta^2(x) \right), \eqno(5.8a)
$$

$$
T_{\zeta_{2,y}}(u_2) = \exp \left( -u_2^2/2 \theta^2(y) \right) \eqno(5.8b)
$$
and are independent. \par

 \ It is no hard to calculate as above the densities of both the independent  variables $ \zeta_{1,x}, \ \zeta_{2,y}.  $ \par

 \ We estimate consequently

$$
\Delta_{n_1,n_2}[f](x,y) = \left| \ B_{n_1,n_2}[f](x,y) - f(x,y) \ \right| = 2^d \times
$$

$$
 \left| {\bf E}(f(x + \zeta_{1,x}/\sqrt{n_1}, y + \zeta_{2,y}/\sqrt{n_2} ) - f(x,y)) \right| \le
2^d {\bf E} \omega[f] \left( \frac{\zeta_{1,x}}{\sqrt{n_1}}, \ \frac{\zeta_{2,y}}{\sqrt{n_2}} \right),
$$
and we deduce after integration twice  by parts analogously to the proof of the inequality (3.5)
$ \  \Delta_{n_1, n_2}[f](x,y) \le 2^d \times $

$$
\int_0^{\infty} \int_0^{\infty} \omega[f] \left( \frac{v_1}{\sqrt{n_1}},
 \frac{v_2}{\sqrt{n_2}}  \right) \ \frac{v_1}{\theta^2(x)} \  \frac{v_2}{\theta^2(y)}
\exp \left( - v_1^2 /2 \theta^2(x) - v_2^2/2\theta^2(y)  \right) \ dv_1 \ dv_2 =
$$

$$
2^d \cdot  \int_0^{\infty}  \int_0^{\infty} \omega[f] \left( \frac{z_1 \cdot \theta(x)}{\sqrt{n_1}}, \
 \frac{z_2 \cdot \theta(y)}{\sqrt{n_2}} \right)
 z_1 \ z_2  \ \exp \left( - \frac{z_1^2}{2} -  \frac{z_2^2}{2} \right) \ dz_1 \ dz_2 =
$$

$$
2^d \ J_{n_1, n_2}[f](x,y). \eqno(5.9)
$$

 \ Thus, we proved the upper bound: $  W_2 \le 2^d.  $ \par

 \ To derive the lower bound, it is sufficient to consider as the capacity  of an example the {\it factorable } function,
namely

$$
f_0(x,y) = g_{t_1}(x) \cdot h(y), \ t_1 \in (0,1), \ h(\cdot) \in C[0,1],
$$
so that

$$
B_{n_1, n_2}[f_0](x,y) =    B_{n_1}[g_{t_1}](x) \cdot B_{n_2}[h](y);
$$

$$
\omega[f_0](\delta_1, \delta_2) = \omega[g_{t_1}](\delta_1) \cdot   \omega[h]( \delta_2),
$$
and finally

$$
J_{n_1, n_2}[f_0](x,y) = J_{n_1}[g_{t_1}](x) \cdot J_{n_2}[h](y).
$$

 \ Let us choose, for example, $ h(y) = \const = 1.  $
 It follows immediately from the second proposition of theorem 3.1

$$
W_2(x,y) \ge \lim_{n_1 \to \infty} \frac{\Delta_{n_1, \infty}[f_0](x,y)}{ J_{n_1,\infty}[f_0](x,y) } =
$$

$$
\lim_{n_1 \to \infty} \frac{\Delta_{n_1}[g_{t_1}](x)}{ J_{n_1}[g_{t_1}](x)} =
 W(x) = 1/\pi, \eqno(5.10)
$$
Q.E.D.\par

\vspace{4mm}

\begin{center}

{\bf  More general approach.} \par

\end{center}

 \ Let $ ||(x,y)|| $ be any non-degenerate norm on the plane $  R^2, $ for instance,

$$
||(x,y)||  = \sqrt{x^2 + y^2}, \ \hspace{4mm} ||(x,y)||  = \max(|x|, |y|), \hspace{4mm} ||(x,y)||  = |x| + |y| \eqno(5.11)
$$
and  so one. \par
 \ In accordance with this definition suppose that there exists non-negative continuous non-decreasing function
$   \gamma(\delta_1, \delta_2) = \gamma[f](\delta_1, \delta_2) $  such that $ \gamma(0+, 0+) = 0 $ and

$$
| f(x_1, y_1) - f(x_2, y_2)| \le \gamma[f](|| (x_1 - x_2, y_1 - y_2) ||).
$$

 \ We deduce as above  $   \Delta_{n_1, n_2}[f](x,y) \le 2^d \times $

$$
\int_0^{\infty} \int_0^{\infty} \gamma[f] \left( \frac{v_1}{\sqrt{n_1}},
 \frac{v_2}{\sqrt{n_2}}  \right) \ \frac{v_1}{\theta^2(x)} \  \frac{v_2}{\theta^2(y)}
\exp \left( - v_1^2 /2 \theta^2(x) - v_2^2/2\theta^2(y)  \right) \ dv_1 \ dv_2. \eqno(5.12)
$$

 \vspace{3mm}

 \section{ Concluding remarks.}

 \vspace{3mm}

 \ It is no hard to generalize the results of the last section into the "more" multivariate case $  d = 3,4,5, \ldots; $
as well as into other methods of approximation, if only they had a probabilistic representation.\par
 \ One can also investigate and improve the rate of convergence of partial derivatives for the multivariate
Bernstein's polynomials. \par

{

\end{document}